\title{Enumerating Regular Objects associated with Suzuki Groups}

\author{Martin Downs and Gareth A. Jones\\
School of Mathematics\\
University of Southampton\\
Southampton SO17  1BJ, U.K.\\
{\tt G.A.Jones@maths.soton.ac.uk}
}

\documentclass[12pt]{article}
\usepackage{color}
\usepackage[latin1]{inputenc}
\usepackage{a4wide}
\usepackage{latexsym}
\usepackage{amsmath}
\usepackage{amssymb}
\usepackage{amsfonts}
\newtheorem{thm}{Theorem}[section]

\newtheorem{prop}[thm]{Proposition}
\date{}
\begin{document}

\maketitle

\begin{abstract}
We use the M\"obius function of the simple Suzuki group $Sz(q)$ to enumerate regular objects such as maps, hypermaps, dessins d'enfants and surface coverings with automorphism groups isomorphic to $Sz(q)$. 
\end{abstract}

\noindent{\bf MSC classification:} 20B25 (primary); % finite automorphism groups
05A15, % exact enumeration problems
05C10, % planar graphs, geometric and topologic aspects of graph theory
05E18, % group actions on combinatorial structures
14H57, % dessins d'enfants theory
20D06 % simple groups, alternating and Lie type
(secondary).

\section{Introduction}

In~\cite{Hal}, Hall introduced the concept of M\"obius inversion in the lattice of subgroups of a group $G$, and used it to find the number $n_{\Gamma}(G)$ of normal subgroups $N$ of a finitely generated group $\Gamma$ with quotient $\Gamma/N$ isomorphic to a given finite group $G$. An important ingredient in this theory is the M\"obius function $\mu_G$ of $G$, a function from the set of subgroups $H$ of $G$ to $\mathbb Z$ defined recursively by
\begin{equation}\label{mu}
\sum_{K\ge H}\mu_G(K)=\delta_{H,G}
\end{equation}
where $\delta_{H,G}$ is the Kronecker delta function, equal to $1$ or $0$ as $H=G$ or $H<G$. Among the groups $G$ for which Hall computed this function were the simple groups $L_2(p)=PSL_2(p)=SL_2(p)/\{\pm I\}$ for primes $p\ge 5$; he then used this to evaluate $n_{
\Gamma}(G)$ for various groups $\Gamma$, such as the free group $F_2$ of rank $2$ and the free products $C_2*C_2*C_2$ and $C_2*C_3\cong PSL_2({\mathbb Z})$ (the modular group). 

In~\cite{DowPhD} the first author extended the calculation of $\mu_G$ to the groups $G=L_2(q)$ and $PGL_2(q)$ for all prime powers $q$; see~\cite{DowJLMS} for $G=L_2(q)$, including a proof for $q=2^e$ and a statement of the corresponding results for odd $q$, together with the calculation of $n_{\Gamma}(G)$ for the modular group $\Gamma$. In~\cite{DJ} the present authors applied these results to enumerate various combinatorial objects, such as regular and orientably regular maps and hypermaps, with automorphism groups isomorphic to $L_2(2^e)$; this is possible since in each of these categories the (isomorphism classes of) regular objects correspond bijectively to the normal subgroups $N$ of some finitely generated group $\Gamma$, with automorphism groups isomorphic to $\Gamma/N$.

In~\cite{Suz}, Suzuki discovered a family of non-abelian finite simple groups $G=Sz(q)={}^2B_2(q)$, where $q=2^e$ for some odd $e>1$, with internal structure similar to that of the groups $L_2(2^e)$. The first author computed the M\"obius function $\mu_G$ for these groups in~\cite{DowSuz}. Our aim here is to apply this to enumerate various regular objects, such as maps, hypermaps, dessins d'enfants and surface coverings, with automorphism groups isomorphic to $G$. This extends work by Silver and the second author~\cite{JS}, where orientably regular maps of type $\{4,5\}$ with automorphism group $G$ were enumerated, and by Hubard and Leemans~\cite{HL}, where regular and chiral maps and polytopes were enumerated. In each case the authors used a restricted form of M\"obius inversion, concentrating mainly on subgroups $H\cong Sz(2^f)$ where $f$ divides $e$. Here we use the full M\"obius function $\mu_G$ (see Table~\ref{table:MobTable} in \S\ref{MobSuz}), allowing a wider range of regular objects to be enumerated; this is done for $G=Sz(q)$ in \S\ref{Enum} and in more detail for $G=Sz(8)$ in \S\ref{Sz8}. A  typical result is that the formula
\[\frac{1}{e}\sum_{f|e}\mu\left(\frac{e}{f}\right)2^f(2^{4f}-2^{3f}-9),\]
where $\mu$ is the classical M\"obius function on $\mathbb N$, counts the normal subgroups of $F_2$ with quotient group $G=Sz(2^e)$, and hence also the orbits of ${\rm Aut}\,G$ on generating pairs for $G$, and the orientably regular hypermaps (or regular dessins d'enfants~\cite{Gro}) with orientation-preserving automorphism group $G$.

\medskip

\noindent{\bf Acknowledgements} The authors are grateful to Dimitri Leemans for some very helpful comments on enumeration with Suzuki groups, and to Nikos Kanakis for help in preparing the TeX file.

\section{M\"obius inversion in groups}

Here we briefly outline Hall's theory of M\"obius inversion in groups~\cite{Hal}. If $\sigma$ and $\phi$ are integer-valued functions defined on isomorphism classes of finite groups, such that
\begin{equation}\label{sigma}
\sigma(G)=\sum_{H\le G}\phi(H)
\end{equation}
for all finite groups $G$, then a simple argument shows that
\begin{equation}\label{phi}
\phi(G)=\sum_{H\le G}\mu_G(H)\sigma(H)
\end{equation}
where $\mu_G$ is the M\"obius function for $G$ defined by equation~(\ref{mu}). We can regard this as inverting equation~(\ref{sigma}) to produce equation~(\ref{phi}); it is an analogue of M\"obius inversion in elementary number theory, where $\sigma$ and $\phi$ are functions ${\mathbb N}\to{\mathbb Z}$, and the corresponding lattice consists of the subgroups $H=n{\mathbb Z}$ of finite index $n$ in $\mathbb Z$, identified with the natural numbers $n\in\mathbb N$ ordered by divisibility.

An important example of such a pair of functions is given by taking $\sigma(G)=|{\rm Hom}(\Gamma, G)|$ and $\phi(G)=|{\rm Epi}(\Gamma, G)|$, where ${\rm Hom}(\Gamma, G)$ and ${\rm Epi}(\Gamma, G)$ are the sets of homomorphisms and epimorphisms $\Gamma\to G$ for some finitely generated group $\Gamma$. These are finite sets, since such a homomorphism is uniquely determined by the images in $G$ of a finite set of generators for $\Gamma$. Since every homomorphism is an epimorphism onto some unique subgroup $H\le G$ we have
\begin{equation}\label{Hom}
|{\rm Hom}(\Gamma,G)|=\sum_{H\le G}|{\rm Epi}(\Gamma,H)|,
\end{equation}
so M\"obius inversion yields
\begin{equation}\label{Epi}
|{\rm Epi}(\Gamma,G)|=\sum_{H\le G}\mu_G(H)|{\rm Hom}(\Gamma,H)|.
\end{equation}
This is often a useful equation, since counting homomorphisms is generally easier than counting epimorphisms. The cost is that one has to evaluate $\mu_G(H)$ for all $H\le G$, but there are two compensations: firstly Hall~\cite{Hal} showed that $\mu_G(H)=0$ if $H$ is not an intersection of maximal subgroups of $G$, so such subgroups $H$ can be deleted from~(\ref{Epi}), and secondly, once $\mu_G$ has been computed it can be applied to many other pairs of functions $\sigma$ and $\phi$. For instance, such a pair could represent the numbers of homomorphisms and epimorphisms $\Gamma\to G$ with certain extra properties, such as being smooth (i.e.~having a torsion-free kernel).

Now let ${\mathcal N}(G)={\mathcal N}_{\Gamma}(G)$ denote the set of normal subgroups $N$ of $\Gamma$ with $\Gamma/N\cong G$, and let
\begin{equation}\label{n(G)defn}
n(G)=n_{\Gamma}(G)=|{\mathcal N}(G)|=|{\mathcal N}_{\Gamma}(G)|.
\end{equation}
These normal subgroups $N$ are the kernels of the epimorphisms $\Gamma\to G$, and two such epimorphisms have the same kernel if and only if they differ by an automorphism of $G$. It follows that $n(G)$ is the number of orbits of ${\rm Aut}\,G$, acting by composition on ${\rm Epi}(\Gamma,G)$. This action is semiregular, since only the identity automorphism of $G$ can fix an epimorphism $\Gamma\to G$; thus all orbits have length $|{\rm Aut}\,G|$ and hence
\begin{equation}\label{n(G)}
n_{\Gamma}(G)=\frac{|{\rm Epi}(\Gamma,G)|}{|{\rm Aut}\,G|}
=\frac{1}{|{\rm Aut}\,G|}\sum_{H\le G}\mu_G(H)|{\rm Hom}(\Gamma,H)|.
\end{equation}

\section{Counting homomorphisms}

In order to apply this method to a specific pair of groups $\Gamma$ and $G$, one needs to be able to count homomorphisms $\Gamma\to H$ for the subgroups $H\le G$. Given a presentation for $\Gamma$ with generators $X_i$ and defining relations $R_j(X_i)=1$, this amounts to counting the solutions $(x_i)$ in $H$ of the equations $R_j(x_i)=1$. For certain groups $\Gamma$, the character table of $H$ gives this information, as illustrated by the following theorem of Frobenius~\cite{Fro}.
\begin{thm}\label{frobchi}
Let ${\mathcal C}_i$ ($i=1, 2, 3$) be conjugacy classes in a finite group $H$. Then the number of solutions of the equation $x_1x_2x_3=1$ in $H$, with $x_i\in{\mathcal C}_i$ for $i=1, 2, 3$, is given by the formula
\begin{equation}\label{trianglechi}
\frac{|{\mathcal C}_1||{\mathcal C}_2||{\mathcal C}_3|}{|H|}\sum_{\chi}\frac{\chi(x_1)\chi(x_2)\chi(x_3)}{\chi(1)}\end{equation}
where $x_i\in{\mathcal C}_i$ and $\chi$ ranges over the irreducible complex characters of $H$.
\end{thm}

If $\Gamma$ is the triangle group
\[\Delta(m_1, m_2, m_3)=\langle X_1, X_2, X_3\mid X_1^{m_1}=X_2^{m_2}=X_3^{m_3}=X_1X_2X_3=1\rangle\]
of type $(m_1, m_2, m_3)$ for some integers $m_i$, then $|{\rm Hom}(\Gamma, H)|$ can be found by summing~(\ref{trianglechi}) over all choices of triples of conjugacy classes ${\mathcal C}_i$ of elements of orders dividing $m_i$. Similarly, the number of smooth homomorphisms $\Gamma\to H$ can be found by restricting the summation to classes of elements of order equal to $m_i$.

When $\Gamma$ is an orientable surface group, that is, the fundamental group
\[\Pi_g=\pi_1{\mathcal S}_g=\langle A_i, B_i\;(i=1,\ldots, g)\mid\prod_{i=1}^g[A_i,B_i]=1\rangle\]
of a compact orientable surface ${\mathcal S}_g$ of genus $g\ge 1$, with $[a,b]$ denoting the commutator $a^{-1}b^{-1}ab$, the following theorem of Frobenius~\cite{Fro} and Mednykh~\cite{Med} is useful (see~\cite{Jon95} for applications):
\begin{thm}\label{FrobMedchi}
In any finite group $H$, the number of solutions $(a_i, b_i)$ of the equation $\prod_{i=1}^g[a_i,b_i]=1$ is given by the formula
\begin{equation}\label{Pi+chi}
|H|^{2g-1}\sum_{\chi}\chi(1)^{2-2g}
\end{equation}
where $\chi$ ranges over the irreducible complex characters of $H$.
\end{thm}

The formula in~(\ref{Pi+chi}) gives $|{\rm Hom}(\Pi_g,H)|$ in terms of the degrees $\chi(1)$ of the irreducible characters of $H$. When $\Gamma$ is a non-orientable surface group
\[\Pi_g^{-}=\langle A_i\;(i=1,\ldots, g)\mid\prod_{i=1}^gA_i^2=1\rangle
\]
of genus $g\ge 1$, the corresponding result of Frobenius and Schur~\cite{FS} is as follows:
\begin{thm}\label{FrobSchchi}
In any finite group $H$, the number of solutions $(a_i)$ of the equation $\prod_{i=1}^ga_i^2=1$ is given by the formula
\begin{equation}\label{Pi-chi}
|H|^{g-1}\sum_{\chi}c_{\chi}^g\chi(1)^{2-g}
\end{equation}
where $\chi$ ranges over the irreducible complex characters of $H$.
\end{thm}
Here $c_{\chi}=|H|^{-1}\sum_{h\in H}\chi(h^2)$ is the Frobenius-Schur indicator of $\chi$, equal to $1, -1$ or $0$ as $\chi$ is respectively the character of a real representation, the real character of a non-real representation, or a non-real character.

\section{Categories and groups}

In some categories $\mathfrak C$, there is a group $\Gamma=\Gamma_{\mathfrak C}$, which we will call the parent group of $\mathfrak C$, such that the set ${\mathcal R}(G)={\mathcal R}_{\mathfrak C}(G)$ of regular objects in $\mathfrak C$ with automorphism group $G$ is in bijective correspondence with the set ${\mathcal N}(G)={\mathcal N}_{\Gamma}(G)$ of normal subgroups of $\Gamma$ with quotient group $G$ (see~\cite{Jon13} for further details). In particular, if $\Gamma$ is finitely generated and $G$ is finite then these two sets have the same finite cardinality
\begin{equation}
r(G)=r_{\mathfrak C}(G):=|{\mathcal R}_{\mathfrak C}(G)|
=n(G)=n_{\Gamma}(G):=|{\mathcal N}_{\Gamma}(G)|,
\end{equation}
so that equation~(\ref{n(G)}) gives
\begin{equation}
r_{\mathfrak C}(G)=\frac{1}{|{\rm Aut}\,G|}\sum_{H\le G}\mu_G(H)|{\rm Hom}(\Gamma,H)|.
\end{equation}

\subsection{Maps, hypermaps and groups}

A map $\mathcal M$ is regular (in the category $\mathfrak M$ of all maps) if its automorphism group $G={\rm Aut}\,{\mathcal M}={\rm Aut}_{\mathfrak M}{\mathcal M}$ acts transitively on vertex-edge-face flags, which can be identified with the faces of the barycentric subdivision $\mathcal B$ of $\mathcal M$. In this case $G$ is generated by automorphisms $r_i\;(i=0,1,2)$ which change (in the only possible way) the $i$-dimensional component of a particular flag, while preserving its $j$-dimensional components for each $j\ne i$. If $\mathcal M$ has type $\{m,n\}$ in the notation of~\cite[Ch.~8]{CM}, meaning that its faces are all $m$-gons and its vertices all have valency $n$, these generators satisfy
\[r_i^2=(r_0r_1)^m=(r_0r_2)^2=(r_1r_2)^n=1.\]
It follows that there is an epimorphism $\theta:\Gamma\to G$, $R_i\mapsto r_i$, where
\begin{equation}\label{GammaM}
\Gamma=\Gamma_{\mathfrak M}
=\langle R_0, R_1, R_2\mid R_i^2=(R_0R_2)^2=1\rangle,
\end{equation}
so $\mathcal M$ determines a normal subgroup $N=\ker\theta$ of $\Gamma$ with $\Gamma/N\cong G$. Conversely, each such normal subgroup determines a regular map $\mathcal M$ with ${\rm Aut}\,{\mathcal M}\cong G$. Two such maps are isomorphic if and only if they correspond to the same normal subgroup, so the set ${\mathcal R}(G)={\mathcal R}_{\mathfrak M}(G)$ of regular maps with automorphism group $G$ is in bijective correspondence with the set ${\mathcal N}(G)={\mathcal N}_{\Gamma}(G)$. If $G$ is finite then the preceding argument gives
\begin{equation}
r_{\mathfrak M}(G)=\frac{1}{|{\rm Aut}\,G|}\sum_{H\le G}\mu_G(H)|{\rm Hom}(\Gamma,H)|.
\end{equation}

This group $\Gamma_{\mathfrak M}$, a free product of its subgroups $\langle R_0, R_2\rangle\cong V_4$ and $\langle R_1\rangle\cong C_2$, can be regarded as the extended triangle group $\Delta[\infty, 2, \infty]$ of type $(\infty, 2, \infty)$, generated by reflections in the sides of a hyperbolic triangle with angles $0, \pi/2, 0$. Other triangle groups play a similar role for related categories. For example, the extended triangle group
\[\Gamma=\Delta[n,2,m]=\langle R_0, R_1, R_2\mid R_i^2=(R_0R_1)^m=(R_0R_2)^2=(R_1R_2)^n=1\rangle\]
is the parent group for maps of all types $\{m',n'\}$ dividing $\{m,n\}$ (meaning that $m'$ divides $m$ and $n'$ divides $n$). For maps of type $\{m,n\}$ one must restrict attention to normal subgroups $N$ of $\Gamma$ such that $R_0R_1$ and $R_1R_2$ have images of orders $m$ and $n$ in $\Gamma/N$.

For the category $\mathfrak H$ of hypermaps, where hyperedges may be incident with any number of hypervertices and hyperfaces, we delete the relation $(R_0R_2)^2=1$ from the presentation~(\ref{GammaM}), giving the group
\[\Gamma_{\mathfrak H}=\Delta[\infty,\infty, \infty]\cong C_2*C_2*C_2.\]
For hypermaps of types dividing $(p,q,r)$, we use the extended  triangle group
\[\Delta[p,q,r]=\langle R_0, R_1, R_2\mid R_i^2=(R_0R_1)^r=(R_0R_2)^q=(R_1R_2)^p=1\rangle.\]

The parent groups for the categories ${\mathfrak M}^+$ and ${\mathfrak H}^+$ of oriented maps and hypermaps are the orientation-preserving subgroups of index $2$ in $\Gamma_{\mathfrak M}$ and $\Gamma_{\mathfrak H}$, generated by the elements $X=R_1R_0$, $Y=R_0R_2$ and $Z=R_2R_1$ satisfying $XYZ=1$. These are the triangle groups
\[\Gamma_{{\mathfrak M}^+}=\Delta(\infty,2,\infty)
=\langle X, Y, Z\mid Y^2=XYZ=1\rangle\cong C_{\infty}*C_2\]
and
\[\Gamma_{{\mathfrak H}^+}=\Delta(\infty,\infty,\infty)
=\langle X, Y, Z\mid XYZ=1\rangle\cong C_{\infty}*C_{\infty}\cong F_2.\]
For oriented maps of types dividing $\{m,n\}$, or oriented hypermaps of types dividing $(p,q,r)$, we use the triangle groups $\Delta(n,2,m)$ and $\Delta(p,q,r)$, restricting attention to torsion-free normal subgroups for maps and hypermaps of these exact types. (See~\cite{JSin} for further background for these categories.)

\subsection{Reflexibility}

The regular objects in the categories ${\mathfrak M}^+$ and ${\mathfrak H}^+$ are often referred to as orientably regular, since they need not be regular  as objects in the larger categories $\mathfrak M$ and $\mathfrak H$. Let $\mathcal H$ be an orientably regular hypermap of type $(p, q, r)$, corresponding to a normal subgroup $N$ of $\Gamma_{{\mathfrak H}^+}=F_2$ with $F_2/N\cong {\rm Aut}_{{\mathfrak H}^+}{\mathcal H}\cong G$ for some group $G$. Then the following are equivalent:
\begin{itemize}
\item $\mathcal H$ is regular in the category $\mathfrak H$ of all hypermaps;
\item $\mathcal H$ has an orientation-reversing automorphism;
\item $N$ is normal in $\Gamma_{\mathfrak H}=C_2*C_2*C_2$;
\item some (and hence each) pair of the canonical generating triple $x, y, z$ for $G$ are inverted by an automorphism of $G$.
\end{itemize}
If these conditions hold we say that $\mathcal H$ is reflexible; otherwise it is chiral, and $\mathcal H$ and its mirror image $\overline{\mathcal H}$ form a chiral pair, isomorphic in $\mathfrak H$ but not in ${\mathfrak H}^+$.

If $\mathcal H$ is reflexible then $\tilde G:={\rm Aut}_{\mathfrak H}{\mathcal H}\cong \Gamma_{\mathfrak H}/N$ is a semidirect product of $G$ by a complement $C_2$ generated by the image $r_i$ of any $R_i\;(i=0, 1, 2)$ in $\tilde G$. The elements of $\tilde G\setminus G$ act by conjugation on the normal subgroup $G$; if one of them induces an inner automorphism then they all do, and we say that $\mathcal H$ is inner reflexible. In this case, $r_i$ induces conjugation by some $g\in G$, so $c:=r_ig$ centralises $G$ and hence $c^2$ is in the centre $Z(G)$ of $G$. If $Z(G)$ is trivial then $\tilde G=G\times C$ where $C=\langle c\rangle\cong C_2$, and there is a non-orientable regular hypermap $\tilde{\mathcal H}={\mathcal H}/C\in{\mathcal R}_{\mathfrak H}(G)$ with orientable double cover $\mathcal H$; this gives a monomorphism ${\mathcal H}\mapsto\tilde{\mathcal H}$, from the inner reflexible maps in ${\mathcal R}_{{\mathfrak H}^+}(G)$ to ${\mathcal R}_{\mathfrak H}(G)$. If, in addition, $G$ has no subgroup of index $2$, then each hypermap in ${\mathcal R}_{\mathfrak H}(G)$ is non-orientable, with an inner reflexible orientable double cover in ${\mathcal R}_{{\mathfrak H}^+}(G)$, so this monomorphism is a bijection. This proves the first part of the following result; the second part is obvious:

\begin{prop}\label{refl}
{\rm (a)} For any finite group $G$ with trivial centre and no subgroup of index $2$, the inner reflexible hypermaps in ${\mathcal R}_{{\mathfrak H}^+}(G)$ are the orientable double covers of the hypermaps in ${\mathcal R}_{\mathfrak H}(G)$; there are $r_{\mathfrak H}(G)$ of them.\newline
{\rm (b)} If, in addition, ${\rm Out}\,G$ has odd order, then every reflexible hypermap in ${\mathcal R}_{{\mathfrak H}^+}(G)$ is inner reflexible, and there are $r_{\mathfrak H}(G)$ of them.
\end{prop}

Every non-abelian finite simple group $G$ satisfies (a), and the Suzuki groups also satisfy (b). The function ${\mathcal H}\mapsto\tilde{\mathcal H}$ preserves types of hypermaps, so the above proposition also applies to maps.

\subsection{Covering spaces}

Under suitable conditions (namely, that $X$ is path connected, locally path connected, and semilocally simply connected~\cite[Ch.~13]{Mun}), the equivalence classes of unbranched coverings $Y\to X$ of a topological space $X$ form a category $\mathfrak C$ in which the connected objects correspond to the conjugacy classes of subgroups of the fundamental group $\Gamma=\pi_1X$; among these, the regular coverings correspond to the normal subgroups $N$ of $\Gamma$, with covering group isomorphic to $\Gamma/N$. If, in addition, $X$ is a compact Hausdorff space, then $\Gamma$ is finitely generated~\cite[p.~500]{Mun}, so one can use the methods described earlier to count regular coverings of $X$ with a given finite covering group. In particular, this applies if $X$ is a compact manifold or orbifold. Indeed, the categories of maps and hypermaps described above can be regarded as obtained in this way from suitable orbifolds $X$, such as a triangle with angles $\pi/p,\pi/q, \pi/r$ for hypermaps of type dividing $(p,q,r)$, or a sphere with three cone-points of orders $p, q, r$ in the oriented case. Similarly, Grothendieck's dessins d'enfants~\cite{GG, Gro} are the finite coverings of a sphere minus three points, so their parent group is its fundamental group $F_2$.

\section{The Suzuki groups}

This section is largely based on Suzuki's description in~\cite{Suz} of the groups named after him; see also~\cite[\S XI.3]{HB} and~\cite[\S 4.2]{Wil}.

\subsection{The definition of the Suzuki Group $G(e)$}
\label{1.1}
Let ${\mathbb F}_q$ be the finite field of $q=2^e$ elements for some odd $e>1$, and let $\theta$ be the automorphism $\alpha\mapsto\alpha^r$ of ${\mathbb F}_q$ where $r=2^{(e+1)/2}$, so that $\theta^2$ is the Frobenius automorphism $\alpha \mapsto \alpha^2$.

For any $\alpha, \beta \in {\mathbb F}_q$ let $(\alpha, \beta)$ denote the $4 \times 4$ matrix
\[
(\alpha,\beta)=
\left(\,\begin{matrix}1&&&\cr \alpha&1&&\cr \alpha^{\theta+1}+\beta&\alpha^{\theta}&\qquad 1&\cr \alpha^{\theta+2}+\alpha \beta+\beta^{\theta}&\beta&\qquad\alpha&\qquad 1\cr\end{matrix}\,\right).
\]
Since
$(\alpha, \beta)(\gamma, \delta) = (\alpha + \gamma, \alpha \gamma^{\theta} + \beta + \delta),$
these matrices $(\alpha, \beta)$ form a group $Q(e)$ of order $q^2$.

The $4\times 4$ diagonal matrices with diagonal entries $\alpha^{1+\theta}$, $\alpha$, $\alpha^{-1}$, $\alpha^{-1-\theta}$ for $\alpha\in{\mathbb F}_q^*:={\mathbb F}_q\setminus\{0\}$ form a cyclic group $A_0(e)$ of order $q-1$. The group $F(e)$ generated by $Q(e)$ and $A_0(e)$ is a semidirect product of a normal subgroup $Q(e)$ by a complement $A_0(e)$, of order $q^2(q-1)$.

Let $\tau$ denote the $4 \times 4$ matrix with entries $1$ on the minor diagonal and $0$ elsewhere. We define $G(e)$ to be the subgroup of $GL_4(q)$ generated by $F(e)$ and $\tau$.
This is the Suzuki group, usually denoted by $Sz(q)$ or ${}^2B_2(q)$. It is, in fact, the subgroup of the symplectic group $Sp_4(q)=B_2(q)$ fixed by a certain automorphism of order $2$.

\subsection{Notation for some subgroups of $G(e)$}
\label{1.2}
Let us fix $e$ and denote
\[G:=G(e), \quad F:=F(e), \quad Q:=Q(e), \quad A_0:=A_0(e).\]
Restricting matrix entries to the subfields ${\mathbb F}_{2^f}$ of ${\mathbb F}_q$ yields a subgroup $G(f)$ of $G$ for each factor $f$ of $e$. If $f$ and $f'$ are factors of $e$ and $f$ divides $f'$ then
\[G(f)\leq G(f'), \quad F(f)\leq F(f'), \quad Q(f)\leq Q(f'), \quad A_0(f)\leq A_0(f').\]

There are cyclic subgroups of $G$ of mutually coprime odd orders
\[2^e\pm r+1=2^e\pm 2^{(e+1)/2}+1,\]
contained in Singer subgroups of $GL_4(q)$ (note that $(2^e+r+1)(2^e-r+1)=q^2+1$ divides $q^4-1$); let us choose a pair of subgroups $A_1, A_2$ of $G$ of these two orders, numbered  according to the rule
\[|A_1|=a_1(e):=2^e+\chi(e)2^{(e+1)/2}+1,\]
\[|A_2|=a_2(e):=2^e-\chi(e)2^{(e+1)/2}+1,\]
where $\chi(e)=1$ or $-1$ as $e\equiv\pm 1$ or $\pm 3$ mod~$(8)$.
If $f$ divides $e$ then $a_i(f)$ divides $a_i(e)$ for each $i=1, 2$, so we define $A_i(f)\;(i=1,2)$ to be the unique subgroup of $A_i$ of order $a_i(f)$. Note that $a_1(e)$ is divisible by $a_1(1)=5$, whereas $a_2(e)$ is not.

Our rule for distinguishing $A_1$ and $A_2$ may seem artificial, and it differs from that used in~\cite{HB, Suz}, where the rule is that $|A_1|>|A_2|$ for all $e$, but it has the advantage, exploited in~\cite{DowSuz}, that if $f$ divides $f'$ then $A_i(f)\leq A_i(f')$ for $i=1, 2$. However, $A_i(f)$ is not necessarily a subgroup of $G(f)$.

\subsection{Properties of some subgroups of $G$.}
\label{1.3}

\begin{enumerate}

\item $G$ has order $q^2(q^2+1)(q-1)$, and is simple if $e > 1$. (The group $G(1)$ is isomorphic to $AGL_1(5)$, of order $20$.)

\item ${\rm Aut}\,G$ is a semidirect product of ${\rm Inn}\,G\cong G$ by a cyclic group of order~$e$ acting as ${\rm Gal}\,{\mathbb F}_q$ on matrix entries, so $|{\rm Aut}\,G|=e|G|$.

\item $G$ acts doubly transitively on an ovoid $\Omega$ of order $q^2+1$ in ${\mathbb P}^3({\mathbb F}_q)$. Its subgroup $F$, the stabiliser of a point $\omega\in\Omega$, acts as a Frobenius group on $\Omega\setminus\{\omega\}$ with kernel $Q$ and complement $A_0$; hence two subgroups of $G$ conjugate to $Q$ intersect trivially, and two subgroups conjugate to $F$ have their intersection conjugate to $A_0$.

\item $Q$ is a Sylow $2$-subgroup of $G$ of order $q^2$ and of exponent $4$. The centre $Z$ of $Q$ consists of the identity and the involutions of $Q$, with $Z\cong Q/Z\cong V_q$.

\item $ZA_0\cong F/Z\cong AGL_1(q)$.

\item The involutions of $G$ are all conjugate, as are the cyclic subgroups of order $4$; however an element of order $4$ is not conjugate to its inverse.

\item All elements of $G$ except those in a conjugate of $Q$ have odd order. Each maximal cyclic subgroup of $G$ of odd order is conjugate to $A_0$, $A_1$ or $A_2$; the intersection of any two of them is trivial.

\end{enumerate}

\subsection{The normalisers of some subgroups of $G$}
\label{1.4}

For $i=1, 2$ the normaliser $B_i$ of $A_i$ is a Frobenius group, with kernel $A_i$ and a complement of order $4$ generated by an element $c_i$ satisfying $c_i^{-1} a c_i = a^{2e}$ for all $a \in A_i$.
For each $f$ dividing $e$ let $B_i(f):=\langle A_i(f), c_i\rangle$, so $|B_i(f)| = 4a_i(f)$. Now let $f > 1$ if $i = 2$. Then $B_i(f)$ is self-normalising, whereas the normaliser of $A_i(f)$ is $B_i$. Similarly, the normaliser $B_0$ of $A_0$ is dihedral. Let $c$ be any involution in $B_0$ and for each $f$ dividing $e$, define $B_0(f):=\langle A_0(f),c\rangle$. If $f > 1$, then $B_0(f)$ is 
self-normalising whereas the normaliser of $A_0(f)$ is $B_0$.

\subsection{Classification of subgroups}
Any subgroup $H\le G$ is a subgroup of some conjugate of $F$ or $B_i\; (i=0, 1, 2)$ or is conjugate to $G(f)$ for some $f$ dividing $e$.

For each $f$ dividing $e$, we have defined the following subgroups of $G$, with the symbols $(f)$ usually omitted when $f=e$:
\[G(f),\; F(f),\; Q(f),\; Z(f),\; B_0(f),\; A_0(f),\; B_1(f),\; B_2(f),\; A_1(f),\; A_2(f).\]
The conjugacy class in $G$ of any of these groups will be denoted by changing the appropriate Roman capital to the corresponding script capital; thus ${\mathcal G}(f), {\mathcal F}, \ldots$ denote the conjugacy classes containing $G(f)$, $F$, and so on.

To show that a set of elements generate $G$, it is sufficient to show that they do not lie in any maximal subgroup of $G$. These are the subgroups in the conjugacy classes
${\mathcal G}(f)$ ($e/f$ prime), $\mathcal F$, ${\mathcal B}_0$, ${\mathcal B}_1$ and ${\mathcal B}_2$.

\section{The M\"obius function of a Suzuki group}
\label{MobSuz}

The first author computed $\mu_G(H)$ for each subgroup $H\le G$ in~\cite{DowSuz}. The non-zero values are given in Table~\ref{MobSuz}; all subgroups $H$ not appearing in Table~\ref{MobSuz} satisfy $\mu_G(H)=0$, and can therefore be ignored in applying equations such as~(\ref{phi}), (\ref{Epi}) and~(\ref{n(G)}). In the third column, $a_i(f)=2^f\pm\chi(f)2^{(f+1)/2}+1$ for $i=1,2$. In the final column, $\mu$ denotes the classical M\"obius function on $\mathbb N$.

\begin{table}[ht]
\centering
\begin{tabular}{| p{2.9cm} | p{2.2cm} | p{3.5cm} | p{2.4cm} |}
\hline
Conjugacy class of $H$ & Number of conjugates & $|H|$ & $\mu_G(H)$ \\
\hline\hline
$\mathcal{G}(f),\; \; \, 1<f\mid e$ & $|G|/|H|$ & $2^{2f}(2^{2f}+1)(2^f-1)$ & $\mu(e/f)$\\
\hline
$\mathcal{F}(f),\; \; 1<f\mid e$ & $|G|/|H|$ & $2^{2f}(2^f-1)$ & $-\mu(e/f)$ \\
\hline
$\mathcal{B}_0(f),\; 1<f\mid e$ & $|G|/|H|$ & $2(2^f-1)$ & $-\mu(e/f)$ \\
\hline
$\mathcal{A}_0(f),\, 1<f\mid e$ & $|G|/2(q-1)$ & $2^f-1$ & $2\frac{(2^e-1)}{(2^f-1)}\mu(e / f)$ \\ [0.6ex]
\hline
$\mathcal{B}_1(f),\; 1<f\mid e$ & $|G|/|H|$ & $4a_1(f)$ & $-\mu(e/f)$ \\
\hline
$\mathcal{B}_2(f),\; 1<f\mid e$ & $|G|/|H|$ & $4a_2(f)$ & $-\mu(e/f)$ \\
\hline
$\mathcal{B}_2(1)$ &  $|G|/2q$ & $4$ & $-2^{e}\mu(e)$ \\
\hline
$\mathcal{B}_0(1)$ & $|G|/q^2$ & $2$ & $-2^{2e-1}\mu(e)$ \\
\hline
$\mathcal I$ & $1$ & $1$ & $|G|\mu(e)$ \\
\hline

\end{tabular}
\caption{Non-zero values of $\mu_G(H)$ for subgroups $H\le G$.}
\label{table:MobTable}
\end{table}

For later use we record here, in Table~\ref{table:|H|k}, the number $|H|_k$ of elements of order $k$ in each of the subgroups $H$ in Table~\ref{table:MobTable}, for $k=2, 4$ and $5$. 

\begin{table}[ht]
\centering
\begin{tabular}{| p{2cm} | p{2.7cm} | p{3.5cm} | p{3cm} |}
\hline
Conjugacy class of $H$  & $|H|_2$ & $|H|_4$ & $|H|_5$ \\
\hline\hline
$\mathcal{G}(f)$ &   $(2^f-1)(2^{2f}+1)$ & $2^f(2^{2f}+1)(2^f-1)$ & $2^{2f}(2^f-1)a_2(f)$\\
\hline
$\mathcal{F}(f)$  & $2^f-1$ & $2^f(2^f-1)$ & $0$ \\
\hline
$\mathcal{B}_0(f)$  & $2^f-1$ & $0$ & $0$ \\
\hline
$\mathcal{A}_0(f)$  & $0$ & $0$ & $0$ \\ 
\hline
$\mathcal{B}_1(f)$  & $a_1(f)$ & $2a_1(f)$ & $4$ \\
\hline
$\mathcal{B}_2(f)$  & $a_2(f)$ & $2a_2(f)$ & $0$ \\
\hline
$\mathcal{B}_2(1)$ & $1$ & $2$ & $0$ \\
\hline
$\mathcal{B}_0(1)$ & $1$ & $0$ & $0$ \\
\hline
$\mathcal I$ & $0$ & $0$ & $0$ \\
\hline

\end{tabular}
\caption{Values of $|H|_k$ for $k=2, 4$ and $5$.}
\label{table:|H|k}
\end{table}

\newpage

\section{Enumerations}
\label{Enum}

We can now use the values of the M\"obius function $\mu_G$ given in Table~\ref{table:MobTable} to enumerate regular objects with automorphism group $G=Sz(q)$ in various categories $\mathfrak C$. Formulae~(\ref{M^+}) and~(\ref{M}), for enumerating maps, have been found in equivalent form by Hubard and Leemans~\cite{HL}, using more direct methods than the general techniques developed here.

\subsection{Orientably regular hypermaps}

If $\mathfrak C$ is the category ${\mathfrak H}^+$ of oriented hypermaps, we take $\Gamma$ to be the free group $F_2$ of rank~$2$. Then $|{\rm Hom}(\Gamma,H)|=|H|^2$ for each subgroup $H\le G$, so
\[
r_{{\mathfrak H}^+}(G)=n_{F_2}(G)=\frac{1}{|{\rm Aut}\,G|}\sum_{H\le G}\mu_G(H)|H|^2.
\]
Now $|{\rm Aut}\,G|=e|G|$, so using the information in Table~\ref{table:MobTable} about the subgroups $H$ of $G$, their orders, numbers of conjugates, and values of $\mu_G(H)$, we obtain, after some routine algebra,
\begin{equation}\label{H^+}
r_{{\mathfrak H}^+}(G)=n_{F_2}(G)=\frac{1}{e}\sum_{f|e}\mu\left(\frac{e}{f}\right)2^f(2^{4f}-2^{3f}-9)
\sim q^5/e.
\end{equation}
(Here we have used the fact that $\Sigma_{f|e}\mu(e/f)=0$ for $e>1$ to eliminate a constant term in the summation.) Formula~(\ref{H^+}) gives the number of orientably regular hypermaps $\mathcal O$ with orientation-preserving automorphism group ${\rm Aut}_{{\mathfrak H}^+}{\mathcal O}\cong G=Sz(q)$, where $q=2^e$ for some odd $e>1$. It also gives the number of regular dessin d'enfants with automorphism group $G$, the number of normal subgroups of the free group $F_2$ with quotient group $G$, and the number of orbits of ${\rm Aut}\,G$ on ordered pairs of generators of $G$. The dominant term in the summation on the right-hand side is the leading term $2^{5f}$ where $f=e$, so simple estimates show that $r_{{\mathfrak H}^+}(G)\sim q^5/e\sim |G|/e$ as $e\to\infty$. (More generally, results of Dixon~\cite{Dix}, Kantor and Lubotzky~\cite{KL}, and Liebeck and Shalev~\cite{LS} on probabilistic generation imply that for all non-abelian finite simple groups, $r_{{\mathfrak H}^+}(G)\sim |G|/|{\rm Out}\,G|$ as $|G|\to\infty$.)

\subsection{Regular hypermaps}

If $\mathfrak C$ is the category $\mathfrak H$ of all hypermaps, then $\Gamma$ is the free product $C_2*C_2*C_2$. Since $G$ cannot be generated by fewer that three involutions, we can restrict attention to smooth homomorphisms and epimorphisms, those that map the three free factors of $\Gamma$ faithfully into $G$. For each $H$ the number of such homomorphisms $\Gamma\to H$ is $|H|_2^3$, where $|H|_2$ is the number of involutions in $H$. The values of $|H|_2$ for the nine conjugacy classes of subgroups $H$ in Table~\ref{table:MobTable} are given in Table~\ref{table:|H|k}, so after some algebra we obtain
\begin{equation}\label{H}
r_{\mathfrak H}(G)=\frac{1}{e}\sum_{f|e}\mu\left(\frac{e}{f}\right)2^f(2^{3f}-2^{2f+1}+2^{f+1}-5)
\sim q^4/e.
\end{equation}
This is the number of regular hypermaps with automorphism group $G$, and also, by Proposition~\ref{refl}, the number of reflexible hypermaps in ${\mathcal R}_{{\mathfrak H}^+}(G)$. Subtracting the formula in equation~(\ref{H}) from that in~(\ref{H^+}) therefore gives the number of chiral hypermaps in  ${\mathcal R}_{{\mathfrak H}^+}(G)$; note that these predominate.

\subsection{Orientably regular maps}

For the category ${\mathfrak M}^+$ of oriented maps we take $\Gamma=C_{\infty}*C_2$. As in the case of hypermaps we may restrict the summation to smooth homomorphisms. There are $|H||H|_2$ such homomorphisms $\Gamma\to H$, so we obtain
\begin{equation}\label{M^+}
r_{{\mathfrak M}^+}(G)=\frac{1}{e}\sum_{f|e}\mu\left(\frac{e}{f}\right)2^f(2^{2f}-2^f-3)
\sim q^3/e.
\end{equation}
(This is equivalent to the formula obtained by Hubard and Leemans in~\cite[Theorem~15]{HL}.) The $k$-valent maps in ${\mathcal R}_{{\mathfrak M}^+}(G)$ correspond to the torsion-free normal subgroups in ${\mathcal N}_{\Gamma}(G)$, where $\Gamma$ is the Hecke group $C_k*C_2$. To count these we consider smooth homomorphisms $\Gamma=C_k*C_2\to H$. There are $|H|_k|H|_2$ of these, so with $k=4$ and $k=5$ for example, Table~\ref{table:|H|k} gives
\begin{equation}\label{M4+}
r_{{\mathfrak M}_4^+}(G)=\frac{1}{e}\sum_{f|e}\mu\left(\frac{e}{f}\right)2^f(2^f-2)\sim q^2/e
\end{equation}
and
\begin{equation}\label{M5+}
r_{{\mathfrak M}_5^+}(G)=\frac{1}{e}\sum_{f|e}\mu\left(\frac{e}{f}\right)(2^f-1)a_2(f)\sim q^2/e,
\end{equation}
where $a_2(f)=2^f-\chi(f)2^{(f+1)/2}+1$. 

A map ${\mathcal M}\in{\mathcal R}_{{\mathfrak M}^+}(G)$ of type $\{n,n\}$, corresponding to a generating triple $(x,y,z)$ for $G$ of type $(n,2,n)$, is self-dual if and only if $G$ has an automorphism transposing $x$ and $z$. Such an automorphism has order $2$ and is therefore inner, induced by conjugation by an involution $i\in G$. This is equivalent to $G$ having a generating triple $(xi,i,x^{-1})$ of type $(4,2,n)$, corresponding to a map ${\mathcal M}^*$ of type $\{n, 4\}$ in ${\mathcal R}_{{\mathfrak M}^+}(G)$. If $\mathcal M$ corresponds to a subgroup $N\in{\mathcal N}_{\Delta}(G)$, where $\Delta:=\Delta(n, 2, n)$, then its median map ${\mathcal M}^{\rm med}$ corresponds to $N$ as an element of ${\mathcal N}_{\Delta^*}(G\times C_2)$, where $\Delta^*:=\Delta(4, 2, n)$ contains $\Delta(n, 2, n)$ with index $2$, and ${\mathcal M}^*={\mathcal M}^{\rm med}/C_2$. The correspondence ${\mathcal M}\mapsto{\mathcal M}^*$ is a bijection, so the number of self-dual maps ${\mathcal M}\in{\mathcal R}_{{\mathfrak M}^+}(G)$ is equal to the number $r_{{\mathfrak M}_4^+}(G)$ of $4$-valent maps ${\mathcal M}^*\in{\mathcal R}_{{\mathfrak M}^+}(G)$, given in~(\ref{M4+}).

\subsection{Regular maps}\label{regmaps}

For the category $\mathfrak M$ of all maps we take $\Gamma=V_4*C_2$. In this case we may restrict attention to homomorphisms which embed the direct factors as subgroups $V$ and $C$, such that  the generator of $C$ commutes with only the identity element of $V$. The only subgroups $H\le G$ containing such subgroups $V$ and $C$ are those conjugate to some $G(f)$, with $V$ and $C$ in the centres of distinct Sylow $2$-subgroups of $H$. Since $G(f)$ has $2^{2f}+1$ Sylow $2$-subgroups, and their centres are elementary abelian group of order $2^f$, one easily obtains
\begin{equation}\label{M}
r_{\mathfrak M}(G)=\frac{1}{e}\sum_{f|e}\mu\left(\frac{e}{f}\right)(2^f-1)(2^f-2)
=\frac{1}{e}\sum_{f|e}\mu\left(\frac{e}{f}\right)2^f(2^f-3)\sim q^2/e.
\end{equation}
As in the case of hypermaps, Proposition~\ref{refl} implies that this is also the number of reflexible maps in ${\mathcal R}_{{\mathfrak M}^+}(G)$, all of them inner reflexible. Subtracting~(\ref{M}) from~(\ref{M^+}) gives the number of chiral maps (see also~\cite[Theorem~16]{HL}), and as before these predominate.

The formulae~(\ref{M}) are the same for the group $G=L_2(q)$, where $q=2^e$. At first this may seem surprising, since $Sz(q)$ is much larger than $L_2(q)$. However, the distribution of involutions in these two groups is very similar, and the above proof can be applied, with only minor changes, to $L_2(q)$.

This proof gives a natural interpretation for the first formula in~(\ref{M}). In $G(f)$, one can assume  by a unique inner automorphism that the generator $R_1$ of $\Gamma$ is sent to $r_1=\tau$, and that $R_0$ and $R_2$ are sent to a pair of  distinct elements $r_0, r_2$ of the form $(0,\beta)$, in the notation of \S5.1.  Then $(2^f-1)(2^f-2)$ is the number of choices for such an ordered pair, the M\"obius inversion over $f$ picks out those triples $(r_i)$ which generate $G$, and division by $e$ counts the orbits of ${\rm Out}\,G$, acting on these triples as ${\rm Gal}\,{\mathbb F}_q\cong C_e$ acts on coefficients $\beta$.

This parametrisation of maps also allows one to determine their types $\{m,n\}$, since $m$ and $n$ are the orders of $r_ir_1$ for $i=0$ and $2$. A matrix $(0,\beta)\tau$ has characteristic polynomial
\begin{equation}
p(\lambda)=\lambda^4+\beta^{\theta}\lambda^3+\beta^2\lambda^2+\beta^{\theta}\lambda+1,
\end{equation}
so its order, as an element $r_ir_1$ of $G$, is the least common multiple of the multiplicative orders of the roots of $p$. Clearly $\lambda=1$ is not a root, so $m$ and $n$ cannot be equal to $2$ or $4$, since elements of $G$ of these orders are unipotent, with all eigenvalues equal to $1$. Thus $m$ and $n$ are both odd. For example, if we take $\beta=1$ then the roots of $p$ are the primitive $5$th roots of $1$, so $r_ir_1$ has order $5$. Specific examples are considered in \S\ref{Sz8}.

None of these regular maps is self-dual. If one were, $G$ would have an automorphism fixing $r_1$ and transposing $r_0$ and $r_2$. This would be induced by an element of $G$ centralising the involutions $r_0r_2$ and $r_1$; however, these lie in distinct Sylow $2$-subgroups of $G$, so their centralisers have trivial intersection.

\subsection{Surface coverings}
\label{surcov}

In order to apply Theorem~\ref{FrobMedchi} to count regular surface coverings with covering group $G$, one needs to know the degrees of the irreducible complex characters of the subgroups $H$ in Table~\ref{table:MobTable}. The irreducible characters of the Suzuki groups $G(f)$ are described in~\cite{Suz} and~\cite[\S XI.5]{HB},
and the degrees for the other subgroups $H$ are easily found; they are given in Table~\ref{table:chardeg}, where $s:=2^f$, $t:=\sqrt{2s}$, $k_i:=(a_i(f)-1)/4$ for $i=1, 2$, and the notation $d^{\{k\}}$ denotes $k$ characters of degree $d$.

\begin{table}[ht]
\centering
\begin{tabular}{| p{1.9cm} | p{1.9cm} | p{6.1cm} | }
\hline
Conjugacy class of $H$ & Conditions on $f$ & Degrees of irreducible characters of $H$  \\
\hline\hline
$\mathcal{G}(f)$ & $1<f\mid e$ & $1, s^2, (s-1)t/2^{\{2\}}, (s^2+1)^{\{(s-2)/2\}}$,\\
&& $(s-1)a_1(f)^{\{k_2\}}, (s-1)a_2(f)^{\{k_1\}}$ \\
\hline
$\mathcal{F}(f)$ & $1<f\mid e$ & $1^{\{s-1\}}, s-1, (s-1)t/2^{\{2\}}$ \\
\hline
$\mathcal{B}_0(f)$ & $1<f\mid e$ & $1, 1, 2^{\{(s-2)/2\}}$ \\
\hline
$\mathcal{A}_0(f)$ & $1<f\mid e$ & $1^{\{s-1\}}$ \\
\hline
$\mathcal{B}_1(f)$ & $1<f\mid e$ & $1,1,1,1, 4^{\{k_1\}}$ \\
\hline
$\mathcal{B}_2(f)$ & $1<f\mid e$ & $1,1,1,1, 4^{\{k_2\}}$ \\
\hline
$\mathcal{B}_2(1)$ & & $1,1,1,1$  \\
\hline
$\mathcal{B}_0(1)$ & & $1,1$  \\
\hline
$\mathcal I$ & & $1$  \\
\hline

\end{tabular}
\caption{Degrees of irreducible characters of subgroups $H\le G$.}
\label{table:chardeg}
\end{table}

With this information, Theorem~\ref{FrobMedchi} gives $|{\rm Hom} (\Gamma, H)|$ for each $H$ in Table~\ref{table:MobTable}, where $\Gamma$ is the fundamental group $\Pi_g$ of an orientable surface ${\mathcal S}_g$ of genus $g$, and then $n_{\Gamma}(G)$ in equation~(\ref{n(G)}) gives $r_g(G)$, the number of regular coverings of ${\mathcal S}_g$ with covering group $G$. The general formulae are very unwieldy, but in \S\ref{Sz8} we will give a simple example.

\section{The smallest simple Suzuki group}
\label{Sz8}

The smallest of the simple Suzuki groups is the group $G=G(3)=Sz(8)$ of order $29120=2^6.5.7.13$. Putting $e=3$ in the enumerative formulae given above, we find that $G$ is the automorphism group of $1054$ regular hypermaps, of which $14$ are maps; all of these are non-orientable. Similarly, it is the orientation-preserving automorphism group of $9534$ orientably regular hypermaps, of which $142$ are maps. By Proposition~\ref{refl}, $1054$ of these orientably regular hypermaps, and $14$ of these orientably regular maps, are reflexible; these are the orientable double covers of the regular hypermaps and maps associated with $G$, so they are all inner reflexible.

\begin{table}[ht]
\centering
\begin{tabular}{| p{1.2cm} | | p{1cm} | p{1cm} | p{1cm} | p{1cm} | p{1cm} |}
\hline
$m\;\setminus\;n$ & $4$ & $5$ & $7$ & $13$ & total \\
\hline\hline
$4$ & $0$ & $4$ & $8$ & $4$ & $16$ \\
\hline
$5$ & $4$ & $4$ & $13$ & $9$ & $30$ \\
\hline
$7$ & $8$ & $13$ & $26$ & $15$ & $62$ \\
\hline
$13$ & $4$ & $9$ & $15$ & $6$ & $34$ \\
\hline
total & $16$ & $30$ & $62$ & $34$ & $142$ \\
\hline

\end{tabular}
\caption{Number of orientably regular maps of type $\{m,n\}$ in ${\mathcal R}_{{\mathfrak M}^+}(Sz(8))$}
\label{table:ORSz8}
\end{table}

\begin{table}[ht]\
\centering
\begin{tabular}{| p{1.2cm} | | p{1cm} | p{1cm} | p{1cm} | p{1cm} | p{1cm} |}
\hline
$m\;\setminus\;n$ & $4$ & $5$ & $7$ & $13$  & total \\
\hline\hline
$4$ & $0$ & $0$ & $0$ & $0$ & $0$ \\
\hline
$5$ & $0$ & $0$ & $1$ & $1$ & $2$ \\
\hline
$7$ & $0$ & $1$ & $2$ & $3$ & $6$ \\
\hline
$13$ & $0$ & $1$ & $3$ & $2$ & $6$ \\
\hline
total & $0$ & $2$ & $6$ & $6$ & $14$ \\
\hline

\end{tabular}
\caption{Number of regular maps of type $\{m,n\}$ in ${\mathcal R}_{\mathfrak M}(Sz(8))$.}
\label{table:RegSz8}
\end{table}

Theorem~\ref{frobchi} and the character table of $G$ in~\cite{ATLAS, Suz} can be used to find how many of these orientably regular maps and hypermaps have a given type. For instance, they show that $G$ contains $2^6.3.7.13.331$ triples $(x,y,z)$ of type $(5,5,5)$ satisfying $xyz=1$; of these, $2^6.3.7.13$ generate the $2^4.7.13$ Sylow $5$-subgroups, while the remaining $2^6.3.7.13.330=66|{\rm Aut}\,G|$ generate $G$, so there are $66$ hypermaps of type $(5,5,5)$ in ${\mathcal R}_{{\mathfrak H}^+}(G)$. Similarly, as shown in~\cite{JS}, there are four maps of type $\{4,5\}$ in ${\mathcal R}_{{\mathfrak M}^+}(G)$, forming two chiral pairs. In fact, repeated use of Theorem~\ref{frobchi} shows that the distribution of types of the orientably regular maps in ${\mathcal R}_{{\mathfrak M}^+}(G)$ is as in Table~\ref{table:ORSz8}, which is symmetric in $m$ and $n$ by the duality of maps. The number of self-dual maps of type $\{n,n\}$ is equal to the number of maps of type $\{n, 4\}$ in this table.

One can use the argument at the end of \S\ref{regmaps} to determine the types $\{m,n\}$ of the $14$ regular maps in ${\mathcal R}_{\mathfrak M}(G)$. Taking $\beta=1$ gives an element $r_ir_1$ of order $5$. Of the six remaining elements $\beta\in{\mathbb F}_8^*$, three have minimal polynomial $t^3+t+1$ over ${\mathbb F}_2$, and three have $t^3+t^2+1$. In the first case $p$ splits into four linear factors, with roots $\beta+1$, $\beta^2$, $\beta^2+\beta$ and $\beta^2+\beta+1$ all of order $7$, so that $r_ir_1$ has order $7$. In the second case, $p$ is irreducible over ${\mathbb F}_8$, and its roots in its splitting field ${\mathbb F}_{2^{12}}$ have order $13$, so $r_ir_1$ has order $13$. By considering the action of ${\rm Gal}\,{\mathbb F}_8\cong C_3$ on distinct ordered pairs of elements $\beta\in{\mathbb F}_8^*$, we find that the number of regular maps of each type $\{m,n\}$ is as in Table~\ref{table:RegSz8}; there are none with $m=4$ or $n=4$ since elements of order $4$ are not conjugate to their inverses. This table also gives the types of the $14$ reflexible maps in ${\mathcal R}_{{\mathfrak M}^+}(G)$, so subtracting its entries from the corresponding entries in Table~\ref{table:ORSz8} gives the number of chiral maps of each type in ${\mathcal R}_{{\mathfrak M}^+}(G)$.

One can also enumerate regular surface coverings with covering group $G$. If we put $f=e=3$ in Table~\ref{table:chardeg}, so that $s=8$ and $t=4$, we find from equation~(\ref{n(G)}) and Theorem~\ref{FrobMedchi} that the number $r_g(G)$ of regular coverings of an orientable surface of genus $g$ with covering group $G$ is
\[\frac{1}{3}\Bigl\{29120^n(1+2.14^{-n}+3.35^{-n}+64^{-n}+3.65^{-n}+91^{-n})\]
\[-448^n(7+7^{-n}+2.14^{-n})-14^n(2+3.2^{-n})+7^{n+1}\]
\[-52^n(4+3.4^{-n})-20^n(4+4^{-n})+8.4^n+2.2^n-1\Bigr\}\]
where $n=2g-2$ is the negative of the Euler characteristic of the surface. When $g=1$ there are no coverings, as one should expect since the fundamental group $\Pi_1$ is abelian, and when $g=2$ there are $286063776$. As $g\to\infty$ we have $r_g(G)\sim |G|^n/3=847974400^{g-1}/3$.

\end{document}